\documentclass[]{amsart}

\usepackage{amsmath,amsthm,amsfonts, bbm, amssymb, hyperref}

\newcommand{\ii}{\mathbbm{i}}
\newcommand{\C}{\mathbb{C}}

\newcommand{\Heis}{\mathcal{H}}
\newcommand{\Hyp}{\mathbb{H}}

\newcommand{\N}{\mathbb{N}}
\newcommand{\R}{\mathbb{R}}
\newcommand{\Z}{\mathbb{Z}}

\newcommand{\CP}{\mathbb{CP}}

\newcommand{\norm}[1]{\left\vert #1 \right \vert}

\renewcommand{\Re}{\text{Re}}
\renewcommand{\Im}{\text{Im}}

\def\[#1\]{\begin{align}#1\end{align}}
\def\(#1\){\begin{align*}#1\end{align*}}

\def\[#1\]{\begin{align}#1\end{align}}

\newtheorem{thm}{Theorem}[section]

\newtheorem{prop}[thm]{Proposition}
\newtheorem{lemma}[thm]{Lemma}
\newtheorem{cor}[thm]{Corollary}
\newtheorem{question}[thm]{Problem}

\theoremstyle{definition}
\newtheorem{defi}[thm]{Definition}

\theoremstyle{definition}
\newtheorem{example}[thm]{Example}
\newtheorem{remark}[thm]{Remark}
\newtheorem{prob}[thm]{Problem}

\numberwithin{equation}{section}

\newcommand{\st}{\;:\;}
\newcommand{\mst}{\;:\;}

\newcommand{\vaisala}{V\"ais\"al\"a}

\begin{document}
\title{Bi-Lipschitz Extension from Boundaries of Certain Hyperbolic Spaces}
\author{Anton Lukyanenko}
\date{\today}
\address{University of Illinois Urbana-Champaign, 1407 W. Green Street, Urbana, IL 61802. Email: Anton@Lukyanenko.net}

\begin{abstract}
Tukia and \vaisala{} showed that every quasi-conformal map of $\R^n$ extends to a quasi-conformal self-map of $\R^{n+1}$. The restriction of the extended map to the upper half-space $\R^n \times \R^+$ is, in fact, bi-Lipschitz with respect to the hyperbolic metric. 

More generally, every homogeneous negatively curved manifold decomposes as $M = N \rtimes \R^+$ where $N$ is a nilpotent group with a metric on which $\R^+$ acts by dilations. We show that under some assumptions on $N$, every quasi-symmetry of $N$ extends to a bi-Lipschitz map of $M$.

The result applies to a wide class of manifolds $M$ including non-compact rank one symmetric spaces and certain manifolds that do not admit co-compact group actions. Although $M$ must be Gromov hyperbolic, its curvature need not be strictly negative.
\end{abstract}

\maketitle

\section{Introduction}
\label{sec:intro}

A classical result of Heintze states that any homogeneous manifold $M$ with negative sectional curvature is diffeomorphic to a Lie group $N \rtimes \R^+$. Here $N$ is a nilpotent group diffeomorphic to $\R^n$ for some $n$, and $\R^+$ acts on $N$ by a family of diffeomorphisms $\alpha_s$. As a negatively curved space, $M$ possesses a boundary at infinity $\partial_\infty M$, composed of geodesic rays up to appropriate equivalence. The abstract boundary $\partial_\infty M$ may be identified with the topological boundary of $M$ as a subset of $(N \times \R) \cup \{\infty\}$; namely, $\partial_\infty M \cong N \cup \{\infty\}$. This identification provides $N$ with a family of parabolic visual metrics; for one of the visual metrics, the diffeomorphism $\alpha_s$ is a homothety with dilation factor $s$.

A map $f: M \rightarrow M$ is a quasi-isometry if there exist $L \geq 1, C \geq 0$ such that for all $x,y \in M$:
\[ \label{eqn:qi}
-C + L^{-1} d_M(x,y) \leq d_M(fx,fy) \leq L d_M(x,y) + C.\]
It is known that every quasi-isometry $f:M\rightarrow M$ extends to a quasi-symmetry $\partial_\infty f:\partial_\infty M\rightarrow \partial_\infty M$. Conversely, every quasi-symmetry of $\partial_\infty M$ extends to a quasi-isometry of $M$, unique up to bounded additive error.

A standard problem is to modify $f$, up to bounded error, so that in (\ref{eqn:qi}) one has $C=0$ or $L=1$. In the first case, $f$ is replaced by a bi-Lipschitz map, and in the second by a rough isometry. In very specialized cases such as Mostow's rigidity theorem, one accomplishes both goals simultaneously: $f$ is bounded distance from an isometry.

In their 1982 paper, Tukia and \vaisala{} examined the bi-Lipschitz approximation problem for real hyperbolic space $\Hyp^n_\R$.
\begin{thm}[Tukia-\vaisala{} \cite{tukiavaisala82}]
\label{thm:tv1}
Let $f: \Hyp^{n+1}_\R \rightarrow \Hyp^{n+1}_\R$ be a quasi-isometry. Then $f$ is bounded distance from a bi-Lipschitz map $F: \Hyp^{n+1}_\R \rightarrow \Hyp^{n+1}_\R$.
\end{thm}

\begin{remark}
The main theorem of \cite{tukiavaisala82} is, in fact, that every quasi-symmetry $g$ of $\R^{n}$ extends to a quasi-symmetry $G$ of $\R^{n+1}$. This follows from Theorem \ref{thm:tv1} by extending $g$ to a quasi-isometry of $\Hyp^{n+1}_\R$ and approximating it by a bi-Lipschitz map.  One then observes that, in the upper half-space model, $\Hyp^{n+1}_\R$ is conformally equivalent to $\R^{n}\times \R^+$ to complete the proof. See Section \ref{sec:tv} for more details.
\end{remark}

Our goal is to generalize Theorem \ref{thm:tv1} to a larger class of hyperbolic spaces. Let $M$ be a Gromov hyperbolic manifold. The following problems are equivalent:
\begin{enumerate}
\item Is every quasi-isometry of $M$ bounded distance from a bi-Lipschitz map?
\item Does every quasi-symmetry of $\partial_\infty M$ extend to a bi-Lipschitz map of $M$?
\end{enumerate}
Following Tukia-\vaisala{}, we focus on version (2) of the problem.

The dyadic tiling of the upper half-space (see Section \ref{sec:tv}) is the essential combinatorial ingredient in the proof of Theorem \ref{thm:tv1}. The tiles are related to each other by isometries of hyperbolic space, and the intersection graph of the tiling captures the geometry of $\Hyp^{n+1}_\R$. We define two classes of spaces (see Section \ref{sec:mss} for precise definitions) that possess properties of $\R^n$ and $\Hyp^{n+1}$ essential to the proof of Theorem \ref{thm:tv1}. 

We say that a metric space $H$ homeomorphic to $\R^n$ posesses a \textit{stacked tiling} if it has an analogue of the integer lattice. Namely, $H$ must admit a discrete co-compact group of isometries $\Gamma$ and a homothety $\alpha$. Furthermore, we require that there exists a fundamental domain $K$ for the action of $\Gamma$ whose rescaling $\alpha K$ is tiled by translates of $K$. The space $H\times \R^+$, endowed with a Riemannian metric, is a \textit{metric similarity space} if both $\Gamma$ and $\alpha$ extend appropriately to isometries of $H^+$. We say that $H$ is the \textit{base} of $H^+$.

One can extend the stacked tiling of the base $H$ to a tiling of $H^+$ by dyadic-type tiles. We show that $H^+$ is Gromov hyperbolic and define the \textit{vertical direction} as the point $[\gamma_\infty]\in \partial_\infty H^+$ represented by quasi-geodesics of the form $\{x\} \times [1,\infty)$, for $x \in H$. We then identify $\partial_\infty H^+ \backslash \{[\gamma_\infty]\}$ with $H \times \{0\}$ as sets, and show that the metric $d_H$ on $H$ is, under the identification, a parabolic visual metric on $\partial_\infty H^+ \backslash \{[\gamma_\infty]\}$. We conclude that any quasi-isometry of $H^+$ preserving the vertical direction extends to a quasi-symmetry of $H$.

We then prove (see Section \ref{sec:main} for the precise statement):
\begin{thm}
\label{thm:main}
Let $H^+=H \times \R^+$ be a metric similarity space of dimension not equal to 4. Then every quasi-symmetry of $H$ extends to a bi-Lipschitz mapping of $H^+$. Equivalently, every quasi-isometry of $H^+$ preserving the vertical direction is bounded distance from a bi-Lipschitz map.
\end{thm}

Theorem \ref{thm:main} applies, in particular, to the homogeneous manifolds $M = N \rtimes \R^+$ where the group structure of the nilpotent group $N$ is defined by polynomials with rational coefficients and the Lie algebra of $N$ is graded (see Definition \ref{ex:graded}).

\begin{thm}
\label{thm:main2}
Let $N$ be a nilpotent group with rational coefficients and graded Lie algebra, with the corresponding Carnot-Carath\`eodory metric. Letting $\R^+$ act on $M$ by homotheties, set $M = N \rtimes \R^+$ and give $M$ a left-invariant Riemannian metric. If the dimension of $M$ is not equal to 4, then any quasi-isometry of $N \rtimes \R^+$ preserving the vertical direction is bounded distance from a bi-Lipschitz map.
\end{thm}

\begin{remark}
Examples of nilpotent groups with rational coefficients and graded Lie algebras include the Heisenberg group and higher-order jet spaces. We may take $M$ to be any non-compact rank one symmetric space of dimension not equal to 4, including complex hyperbolic spaces $\Hyp^n_\C$, $n \neq 2$. Note that for $\Hyp^n_\C$, the condition of preserving the vertical direction is not necessary because the isometry group acts transitively on $\partial_\infty \Hyp^n_\C$.
\end{remark}

Another generalization of Theorem \ref{thm:tv} was given recently by X. Xie. 
\begin{thm}[Xie \cite{xie2009}]\label{thm:xie1}
Let $\widetilde  M$ be the universal cover of a compact negatively curved manifold of dimension not equal to 4. Then every quasi-isometry of $\widetilde  M$ is bounded distance from a bi-Lipschitz map.
\end{thm}

Both Theorem \ref{thm:main} and Theorem \ref{thm:xie1} apply to spaces with rich isometry groups, such as $\Hyp^n_\C$. Apart from the core components, both the assumptions and proofs of the two theorems are quite different, as we discuss in Section \ref{sec:xie}.

The structure of the paper is as follows. We first review the theory of Gromov hyperbolic spaces and their boundaries in Section \ref{sec:defs}. In Section \ref{section2} we summarize the proof of Theorem \ref{thm:tv}, define and study metric similarity spaces, and prove Theorem \ref{thm:main}. In Section \ref{sec:discussion}, we provide examples of metric similarity spaces and compare our result to Theorem \ref{thm:xie1}, concluding with some remarks and open questions.

I would like to thank my PhD adviser Jeremy Tyson for his help and guidance on this project, as well as Richard Schwartz and John Mackay for interesting conversations on the subject. The paper also benifited highly from the 2011 summer school on filling invariants and asymptotic structures organized by David Fisher.

\section{Gromov Hyperbolic Geometry}
\label{sec:defs}

\subsection{Basics}
\label{sec:gromov}

We briefly review the theory of Gromov hyperbolic spaces. For a more complete exposition, see \cite{bonk-schramm, ABC1991, ghysdelaharpeBook} and references therein.

For $C>0$, a \textit{$C$-rough isometry} between two metric spaces $(X,d_X)$ and $(Y,d_Y)$ is a map $f:X \rightarrow Y$ such that
\[-C + d_X(x,y) \leq d_Y(fx, fy) \leq  d_X(x,y) + C\]
for all $x,y \in X$. Note that a rough isometry need not be continuous.

A ($C$-rough) \textit{geodesic} in a metric space $X$ is an ($C$-rough) isometry $f: I \rightarrow X$, where $I$ is an interval. The space $X$ is  ($C$-roughly) \textit{geodesic} if any two points are joined by a  ($C$-rough) geodesic. If $X$ is $C$-roughly geodesic and $C$ is clear from context, we denote any $C$-rough geodesic joining points $x, x' \in X$ by $[x,x']$, even though it may not be unique.

Let $(X,d)$ be a metric space and $x, x', y \in X$. The Gromov product $(x\vert x')_y$ measures the failure of the triangle inequality to be an equality.
\[ (x \vert x')_y = \frac{1}{2} \left( d(x , y) + d(y, x') - d(x, x')\right).\]

\begin{defi}
Let $\delta\geq 0$. A roughly geodesic metric space $X$ is $\delta$-hyperbolic if for all $x,x',x'',y \in X$ one has 
\[
\label{fla:hyperbolic}
(x \vert x'')_y \geq \min\{ (x \vert x')_y, (x' \vert x'')_y\} - \delta.\]
$X$ is called \textit{Gromov hyperbolic} or just \textit{hyperbolic} if it is $\delta$-hyperbolic for some $\delta$. 
\end{defi}

In the context of geodesic metric spaces, one uses the following as the definition of hyperbolicity.

\begin{defi}
A geodesic space $X$ has \textit{$\delta$-thin triangles} if for all points $x,y,z \in X$ and any choice of geodesics $[x,y], [y,z], [x,z]$, each side of the triangle is in the $\delta$-neighborhood of the other two sides:
\[[x,z] \subset N_\delta([x,y] \cup [y,z]).\]
$X$ is said to have thin triangles if its triangles are $\delta$-thin for some $\delta$.
\end{defi}

\begin{lemma}[see \cite{ABC1991}]
\label{lemma:hyperthin}
Let $X$ be a geodesic space. Then $X$ has thin triangles if and only if it is hyperbolic.
\end{lemma}

Bonk-Schramm linked the study of roughly geodesic hyperbolic spaces to that of geodesic hyperbolic spaces with the following theorem.

\begin{thm}[Bonk-Schramm \cite{bonk-schramm}]
Let $X$ be a roughly geodesic $\delta$-hyperbolic space. Then $X$ embeds isometrically into a geodesic $\delta$-hyperbolic space.
\end{thm}

\begin{cor}
\label{cor:roughtriangles}
Let $X$ be a $C$-roughly geodesic $\delta$-hyperbolic space. There exists $\delta'>0$ so that all $C$-roughly geodesic triangles of $X$ are $\delta'$-thin.
\begin{proof}
It is easy to show that there exists a $D>0$ depending only on $C$ and $\delta$ so that in every $\delta$-hyperbolic geodesic metric space, any $C$-rough geodesic is at most distance $D$ from a geodesic with the same endpoints.

Consider an embedding $X \hookrightarrow Y$ of $X$ into a geodesic $\delta$-hyperbolic space. A $C$-roughly geodesic triangle in $X$ can then be approximated by a geodesic triangle in $Y$ with the same vertices.  By Lemma \ref{lemma:hyperthin}, the approximating triangle is thin, so the original $C$-roughly geodesic triangle must also be thin.
\end{proof}
\end{cor}

\subsection{Boundaries of Hyperbolic Spaces}
\begin{defi}
Let $X$ be a roughly geodesic hyperbolic space. The \textit{boundary at infinity} $\partial_\infty X$ is defined as the set of all rougly geodesic rays in $X$, considered equivalent if they are a bounded Hausdorff distance from each other.
\end{defi}

\begin{remark}
There are three additional standard ways of defining the boundary at infinity:
\begin{enumerate}
\item Geodesic rays in $X$, equivalent if they are at bounded Hausdorff distance from each other.
\item Quasi-geodesic rays in $X$, with the same equivalence.
\item Sequences of points in $X$ whose Gromov product converges to $\infty$. Two sequences are equivalent if their Gromov product converges to $\infty$.
\end{enumerate}

For equivalence of the three definitions in geodesic spaces, see \cite{ABC1991}. Equivalence with our definition follows from Section 5 of \cite{bonk-schramm}.
\end{remark}

\begin{defi}
Let $X$ be a roughly geodesic hyperbolic space, $p\in X$, and $a>0$. The \textit{($a,p$)-visual quasimetric} on $\partial_\infty X$ is defined by
\[
\label{fla:previsual}
d_{p,a}\left([\xi],[\zeta]\right) =\lim_{t\rightarrow \infty} a^{-(\xi(t) \vert \zeta(t))_p}.\]
Any metric bi-Lipschitz to $d_{p,a}$ for some $p$ and $a$ is called an \textit{($a$-)visual metric}.
\end{defi}
It is well-known that the limit in (\ref{fla:previsual}) exists. There exists an $a_0>1$ so that for any $a\geq a_0$ the function $d_{p,a}$ is bi-Lipschitz to a metric (see \cite{ghysdelaharpeBook} Chapter 7). For CAT(-1) spaces, $d_{p,a}$ is itself a metric. For different basepoints $p$, the corresponding functions $d_{p,a}$ are bi-Lipschitz equivalent, so that the family of visual metrics does not depend on $p$.

\begin{defi}
Let $X$ be a roughly geodesic hyperbolic space, $\eta$ a roughly geodesic ray, and $a>0$. The \textit{($[\eta],a$)-parabolic visual quasimetric} is defined by
\[
\label{fla:horo}
d_{[\eta],a}\left([\xi] \vert [\zeta]\right) = \lim_{t\rightarrow \infty} a^{t-(\xi(t), \zeta(t))_{\eta(t)}}.\]
Any metric on $\partial_\infty X \backslash \{[\eta]\}$ bi-Lipschitz to $d_{[\eta], a}$ for some $a$ is called an \textit{($a$)-parabolic visual metric}.
\end{defi}

As for visual metrics, there exists an $a_0$ so that for any $a>a_0$ the function $d_{[\eta],a}$ is bi-Lipschitz to a metric (see \cite{ghysdelaharpeBook} Chapter 8). 

\begin{remark}
In general, there is no way to identify punctured boundaries corresponding to different punctures $[\eta]\in\partial_\infty X$. While visual metrics based inside $X$ have bounded diameter (since the Gromov product is non-negative), parabolic visual metrics are unbounded, with the puncture $[\eta]$ positioned infinitely far from the other points. Shanmugalingam-Xie related parabolic visual metrics to visual metrics by means of metric inversions in \cite{shanmugalingam-xie2011}.
\end{remark}

Because the Gromov product is bounded below by $0$, visual metrics are bounded. Bonk-Schramm show that every complete bounded metric space arises as the visual boundary of some hyperbolic space (see Section \ref{sec:interior}). We are primarily interested in the case of manifolds with large isometry groups. In this case, one has more control on the boundary.

Heintze showed in \cite{heintze74} that every homogeneous manifold $M$ with negative curvature is diffeomorphic to $N \rtimes \R^+$ with a left-invariant metric, where $N$ is a nilpotent Lie group. There is a left-invariant metric $d_N$ on $N$ so that in the semi-direct product structure $\R^+$ acts on $N$ by homotheties. It is known (cf. Theorem \ref{thm:boundary}) that $\partial_\infty M$ can be identified with the one-point compactification of $N$ and that under this identification $d_N$ is a parabolic visual metric.  As a nilpotent Lie group, $N$ is diffeomorphic to $\R^n$ for some $n$, so $\partial_\infty M$ is a sphere. 

Of particular interest is the case when $M$ is a non-compact rank one symmetric space. For real hyperbolic space $\Hyp^{n+1}_\R$, the round metric on $S^n$ serves as a visual metric on $\partial_\infty \Hyp^{n+1}_\R$. The Euclidean metric on $\R^n$ is a parabolic visual metric on the punctured boundary of $\Hyp^{n+1}_\R$. For complex hyperbolic space $\Hyp^{n+1}_\C$, these metrics are no longer visual. Instead, one gives $S^{2n+1}$ the sub-Riemannian metric induced by the inclusion $S^{2n+1} \hookrightarrow \C^{2n+2}$ to obtain a visual metric on $\partial_\infty \Hyp^{n+1}_\C$. Likewise, $\R^{2n+1}$ is given the stucture of the Heisenberg group and the corresponding Carnot-Carath\`eodory metric to obtain a parabolic visual metric on the punctured boundary. See Section \ref{sec:rankone} for more details.

\subsection{Quasi-Isometries and Quasi-Symmetries}
\label{sec:qiqs}
\begin{defi}
Let $(X,d_X)$ and $(Y,d_Y)$ be metric spaces, $C \geq 0$, and $L \geq 1$. An \textit{$(L,C)$-quasi-isometry from $X$ to $Y$} is a map $f: X \rightarrow Y$ such that for all $x,y \in X$
\[\label{fla:qi}
-C+L^{-1}d_X(x,y) \leq d_Y(fx,fy) \leq L d_X(x,y) + C
\]
and furthermore $Y$ is covered by an $L$-neighborhood of $f(X)$. If the latter condition is not satisfied, we say $f$ is a quasi-isometry \textit{into} $Y$. Note that in both cases $f$ need not be continuous, injective or surjective.
\end{defi}

\begin{defi}
Let $f:X\rightarrow Y$ be a quasi-isometry between $\delta$-hyperbolic metric spaces. If $Y$ is geodesic and proper, then the Arzela-Ascoli theorem implies that the image of every geodesic ray under a quasi-isometry is uniformly close to a geodesic ray. One may then define the induced map $\partial_\infty f$ on the boundary as
$\partial_\infty f([\gamma]) := [\gamma']$, where $\gamma'$ is some geodesic ray close to $f(\gamma)$. If $Y$ is not proper, the Arzela-Ascoli Theorem does not apply, and the construction of $\partial_\infty f$ is more involved; see \cite{bonk-schramm}. 
\end{defi}

\begin{defi}
Let $(X,d_X)$ and $(Y,d_Y)$ be metric spaces and $\eta$ an increasing bijection $\eta: [0,\infty) \rightarrow [0,\infty)$. 
A surjective map $g: X \rightarrow Y$ is an \textit{$\eta$-quasi-symmetry} if for all $x,y,z \in X$ we have
\[ 
\frac{d(gx, gy)}{d(gx,gz)} \leq \eta \left( \frac{d(x,y)}{d(x,z)} \right).
\]
\end{defi}

\begin{remark}
For $X$ and $Y$ equal to Euclidean space or the Heisenberg group, quasi-symmetry is equivalent to quasi-conformality.
\end{remark}

Note that for a hyperbolic space $X$, any two visual metrics on $\partial_\infty X$ are quasi-symmetric to each other. However, the class of visual metrics is not closed under quasi-symmetries. 
\begin{defi}
Let $X$ be a hyperbolic space. A metric $d$ on $\partial_\infty X$ is in the \textit{conformal gauge} if there exists a visual metric $d_{\text{vis}}$ on $\partial_\infty X$ such that the identity map $id: (\partial_\infty X, d) \rightarrow (\partial_\infty X, d_{\text{vis}})$ is a quasi-symmetry. The conformal gauge on the punctured boundary $\partial_\infty X \backslash \{ [\xi]\}$ is, likewise, the class of metrics quasi-symmetrically equivalent to a parabolic visual metric. Note that visual metrics are not quasi-symmetric to parabolic visual metrics, since quasi-symmetries preserve boundedness. They may, however, be quasi-conformally equivalent.
\end{defi}

The following is a classical result for geodesic spaces (see \cite{margulis, bourdon-shadows}). For roughly geodesic spaces, see Theorem 6.5 of \cite{bonk-schramm}.
\begin{thm}
\label{thm:qiqs}
Let $X$ be a hyperbolic space and $f: X \rightarrow X$ a quasi-isometry. Then $\partial_\infty f$ is a quasi-symmetry with respect to any metric in the conformal gauge of $\partial_\infty X$. If some point $[\xi]\in \partial_\infty f$ is fixed by $\partial_\infty f$, then $\partial_\infty f$ is a quasi-symmetry with respect to any metric in the conformal gauge of $\partial_\infty X \backslash \{[\xi]\}$.
\begin{proof}[Sketch of proof]
One shows that $\partial_\infty f$ is a quasi-symmetry (in the appropriate sense) for the quasimetrics $d_{p,a}$ and $d_{[\xi],a}$. It is then a quasi-symmetry with respect to any metric in the conformal guage.
\end{proof}
\end{thm}

For visual metrics, Bonk-Schramm provided a partial converse to Theorem \ref{thm:qiqs} in  \cite{bonk-schramm} by showing that every quasi-symmetry of $\partial_\infty X$ is induced by a quasi-isometry of $X$ itself. Our main theorem states that under further assumptions on $X$, every quasi-symmetry of $\partial_\infty X \backslash {[\xi]}$ is in fact induced by a bi-Lipschitz map on $X$.

\subsection{Recovering the interior}
\label{sec:interior}
Does every metric space arise as the boundary of some hyperbolic space?

\begin{defi}
\label{defi:conx}
Let $(X,d_X)$ be a metric space with diameter $0 < D \leq \infty$. Define a new metric space $\text{Con}(X)$ with underlying set $\text{Con}(X) = X \times (0, D]$ (or $(0, \infty)$ if $D=\infty$)  and metric:
\[\label{fla:dCon}  &d_{\text{Con}(X)}\left( (x,t), (x', t') \right) = 2 \log \left( \frac{ d_X(x,x') + \max\{t, t'\}}{\sqrt{t t'}} \right). \]
\end{defi}

The following theorem combines several results in Sections 7 and 8 of \cite{bonk-schramm}:
\begin{thm}[Bonk-Schramm \cite{bonk-schramm}] 
\label{thm:bonk-schramm}
Suppose $(X, d_X)$ is a complete metric space of diameter $D \leq \infty$. Then $\text{Con}(X)$ is hyperbolic and roughly geodesic. If $D<\infty$, then under a natural identification of $X$ with $\partial_\infty \text{Con}(X)$, the metric $d_X$ is $e$-visual.
\end{thm}
We now suppose that $X$ is unbounded and generalize the identification in Theorem \ref{thm:bonk-schramm}. The following are geodesic rays in Con$(X)$:
\[
\label{fla:gammax}
&\gamma_x(s)  = (x, e^{-s}) &\gamma_{\infty,x_0}(s) = (x_0,e^s),\]
where $s\in [0, \infty)$, and $x, x_0 \in X$. For different choices of $x_0$, the rays $\gamma_{\infty, x_0}$ are bounded distance from each other. When the specific $x_0$ is understood or irrelevant, we will simply write $\gamma_\infty$. Given a roughly geodesic ray $\gamma$, denote its equivalence class in the boundary by $[\gamma]\in \partial_\infty X$. We will refer to $[\gamma_\infty]$ as the \textit{vertical direction}. 

Consider now $\text{Con}(X) = X \times \R^+$ as a subset of $(X\times \R)\cup \{\infty\}$. Identifying a geodesic $\gamma_x$ with its endpoint on the boundary gives the mapping
\[
\label{fla:bdbijection}
\iota: \overline X= X \cup \{ \infty \} &\rightarrow \partial_\infty \text{Con}(X)\\ x&\mapsto [\gamma_x].\]
It is clear that $\iota$ is well-defined and injective. We next show that for some $X$ it is a bijection.

\begin{prop}
\label{prop:boundary}
Let $(X,d_X)$ be an unbounded complete metric space.  Suppose there is a homothety $\alpha$ dilating $X$ by a factor $a>1$ so that for some compact $K\subset X$, $\cup_{n \in \Z} \alpha^n K = X$. Then $\iota$ is a bijection.
\begin{proof}
The metric $d_X$ restricts to a metric on $K$, which we continue to denote $d_X$.  We also have a natural isometric embedding $\text{Con}(K) \hookrightarrow \text{Con}(X)$. Abusing notation, extend $\alpha$ to an isometry of $\text{Con}(X)$ by taking $\alpha(x,t)=(\alpha(x), a t)$. Then, $\text{Con}( \alpha^n K) \subset \text{Con}(X)$ for each $n$, isometric to $\text{Con}(K)$ by $\alpha^n$. Each $(\alpha^n K, d_X)$ is complete and bounded. By Theorem \ref{thm:bonk-schramm} and Corollary \ref{cor:roughtriangles}, there exist $C,C'\geq C$ so that $\text{Con}(X)$ is $C$-rough geodesic and $C$-rough triangles in $\text{Con}(X)$ are $C'$-thin. 

Suppose that $\iota$ is not surjective. Then there exists a rough geodesic ray $\xi$ that is not equivalent to $\gamma_x$ for any $x\in \overline X$. Increasing $C$ and $C'$ if necessary, we may assume that $\xi$ is $C$-rough and has $\xi(0)=(x_0,1)$ for some $x_0 \in X$. Write $\xi(s) = (x(s), t(s))$ for $x(s) \in X$, $t(s)\in \R^+$.

We first show that $t(s)$ is bounded above. Suppose $t(s)$ is unbounded. Let $\gamma_\infty=\gamma_{\infty,x_0}$. Since $[\xi] \neq [\gamma_\infty]$, there exists $s_0>0$ so that $d(\gamma_\infty(s_0), \xi)> C'$. If $t(s)$ is unbounded, we can select $s_1$ so that 
\[\label{fla:s1}
\log(t(s_1)) > s_0 + C + C'.\]
Set $s_2=\log t(s_1)$ so that we have $\gamma_\infty(s_2) = (x_0, t(s_1))$. Fix a $C$-rough geodesic $[\xi(s_1),\gamma_\infty(s_2)]$.

Consider now the triangle with vertices $(x_0,1), \xi(s_1)$, and $\gamma_\infty(s_2)$ and edges $\xi, \gamma_\infty$, and $[\xi(s_1),\gamma_\infty(s_2)]$. Since this triangle is $C'$-thin, by choice of $s_0$ we have 
\[
\label{fla:thin}
d(\gamma_\infty(s_0), [\xi(s_1),\gamma_\infty(s_2)])<C'.\]
We compute some distances using (\ref{fla:dCon}):
\(d(\xi(s_1), \gamma_\infty(s_2)) &= 2 \log\left( \frac{d( x(s_1), x_0) + t(s_1)}{t(s_1)}\right)\\
d( \xi(s_1), \gamma_\infty(s_0)) &= 2 \log\left( \frac{ d(x(s_1), x_0) + t(s_1)}{\sqrt{t(s_1) e^{s_0}}}\right)\\ 
&= 2 \log\left( \frac{d( x(s_1), x_0) + t(s_1)}{t(s_1)}\right) + \log(t(s_1)/e^{s_0})
\)
Applying (\ref{fla:s1}), we conclude
\(d( \xi(s_1), \gamma_\infty(s_0)) > d(\xi(s_1), \gamma_\infty(s_2)) + C + C'.\)
Thus, the $C$-rough geodesic between $\xi(s_1), \gamma_\infty(s_2)$ cannot come within $C'$ of the point $\gamma_\infty(s_0)$, contradicting \ref{fla:thin}. 

We thus have that $t(s)$ is bounded above. Now, by Theorem \ref{thm:bonk-schramm}, we know $\xi$ is not contained in $\alpha^n \text{Con}(K)$ for any $n$. Thus, the rough geodesic ray $\xi$ must travel arbitrarily far along $X$ while maintaining a bounded $t$ coordinate. Such a path is also inefficient, as can be shown with a similar calculation as above. Thus, $\xi$ cannot exist and $\iota$ is a bijection.

\end{proof}
\end{prop}

\begin{cor}
\label{cor:horometric}
Let the notation be as in Proposition \ref{prop:boundary} and identify $\overline X$ with $\partial_\infty \text{Con}(X)$ by means of $\iota$. Assume further that $\alpha$ has a fixed point $x_0\in X$ and that $x_0$ is an interior point of $K$. Then the metric $d_X$ is a parabolic visual metric on $\partial_\infty \text{Con}(X)\backslash \{[\gamma_\infty]\}$.
\begin{proof}
Under $\iota$, $K \subset X$ is identified with $\partial_\infty\text{Con}(K) \subset \partial_\infty\text{Con(X)}$. 
By Theorem \ref{thm:bonk-schramm}, the metric $d_X$ restricted to $K \subset X$ is a visual metric with base $e$. We may take $o=(x_0, 1)$ as the basepoint for the Gromov products in (\ref{fla:previsual}).

Let $x,y \in X$. For large $n$ we have $\alpha^{-n} x, \alpha^{-n} y \in K$ since the distances to $x_0 \in K$ can be made arbitrarily small. We compute for large $n\in \N$, denoting bi-Lipschitz equivalence independent of $n$ by $\approx$.
\(d_X(x, y) &= a^n d_X(\alpha^{-n}x, \alpha^{-n}y) &\text{($\alpha$ is an $a$-homothety)}\\
						 &\approx \lim_{s\rightarrow \infty} a^n e^{-( \alpha^{-n}\gamma_x(s) \vert \alpha^{-n}\gamma_y(s))_o} &\text{(definition of visual metric)}\\
						 &= \lim_{s\rightarrow \infty} a^n e^{-(\gamma_x(s) \vert \gamma_y(s))_{ \alpha^n o}} &\text{(applying the isometry $ \alpha^n$)}\\
						 &= \lim_{s\rightarrow \infty} e^{n \log a - (\gamma_x(s) \vert \gamma_y(s))_{\alpha^n o}}&\text{(change of base)}
\)
Note now that $\alpha^n o = \gamma_\infty(n \log a)$. We thus have
\( d_X(x,y) \approx \lim_{s \rightarrow \infty} e^{n \log a - (\gamma_x(s) \vert \gamma_y(s))_{\gamma_\infty(n \log a)}}.
\)						 
We may take $n = \lfloor s \rfloor / \log a$, so that
\( d_X(x,y) \approx \lim_{s\rightarrow \infty} e^{s - (\gamma_x(s) \vert \gamma_y(s))_{\gamma_\infty(s)}}.
\)
Comparing to (\ref{fla:horo}), $d_X$ is a parabolic visual metric.

\end{proof}
\end{cor}

\section{Metric Similarity Spaces}
\label{section2}
\subsection{Motivating Example: Real Hyperbolic Space}
\label{sec:tv}
We now sketch a proof of the theorem of Tukia and \vaisala{} that motivated this paper. 

\begin{thm}[\cite{tukiavaisala82}]
\label{thm:tv}
Let $f: \R^n \rightarrow \R^n$ be quasi-conformal. Then there exists a quasi-conformal $F: \R^{n+1} \rightarrow \R^{n+1}$ such that the restriction to $\R^{n}\times \{0\}$ agrees with $f$. Furthermore, $F$ restricted to the upper half-space is bi-Lipschitz with respect to the hyperbolic metric.
\begin{proof}[Proof sketch]
Quasi-conformal maps on $\R^n$ are quasi-symmetric. Fix $\eta$ so that $f$ is $\eta$-quasi-symmetric. Denote $H^+=\R^n \times \R^+$, with the hyperbolic metric (see Section \ref{sec:rankone}). Let $K=[0,1]^n$, $Q_0 = K \times [1,2]$. Then the elements of  
$\mathcal Q = \{ 2^{a}(Q_0 + b) \st a \in \Z, b \in \Z^n\}$ are the \textit{tiles} of the dyadic decomposition of $H^+$. With respect to the hyperbolic metric, the tiles are isometric to each other by maps that preserve the local combinatorics of the tiling.

Given a homeomorphism $g: \R^n \rightarrow \R^n$, define the \textit{basic extension} $\widetilde  g: H^+ \rightarrow H^+$ by
\(\widetilde  g (x, t) = \left(g(x), \max_{\norm{x-y} \leq t} \norm{g(x)-g(y)}\right).\)
The extension $\widetilde  g$ is a homeomorphism that extends continuously to $\partial H^+$, with values on $\partial H^+$ given by $g$. The extension $\widetilde  f$ of $f$ is not necessarily bi-Lipschitz on small scales, but this can be repaired using a combination of piecewise-linear (PL) approximation and Sullivan's extension theorem, as we now describe. 

Recall that, except in dimension 4, any homeomorphism of Euclidean space may be approximated by a PL homeomorphism. Because $H^+$ is Riemannian, it is locally bi-Lipschitz to Euclidean space, and Euclidean methods such as PL approximation may be applied locally. 

We would like to apply PL approximation to obtain a bi-Lipschitz approximation $F$ of $f$ on each tile $Q \in \mathcal Q$ . However, the approximations on adjacent tiles need not agree. To avoid this, let $\mathcal Q = \mathcal Q_1 \sqcup \ldots \sqcup \mathcal Q_N$ be a coloring of $\mathcal Q$. On tiles in $\mathcal Q_1$, define $F$ using a PL approximation of $\widetilde  f$. 

\begin{remark}
In the case $n=4$, one can first use PL approximation on $\R^n$ and then lift to $H^+$ while preserving the bi-Lipschiz property. Note that there do exist homeomorphisms of $\R^4$ to itself that cannot be piecewise-linearly approximated.
\end{remark}

The next step is to interpolate the approximations. This is done using Sullivan's bi-Lipschitz extension result (see Lemma \ref{lemma:BL-extension}). For each tile in $\mathcal Q_2$, Sullivan extension provides values of $F$ that agree with the values given on the tiles $\mathcal Q_1$. Once the interpolation is complete on $\mathcal Q_1$, continue inductively on the remaining colors until $F$ is defined on all of $H^+$.

Every step above used infinitely many bi-Lipschitz approximations. To guarantee that the bi-Lipschitz constant stays bounded, fix $p_0 \in H^+$ and let
\(\widetilde{\mathcal F_\eta} = \{ \widetilde  g: H^+ \rightarrow H^+ \mst \widetilde  g(p_0) \in Q_0 \text{ and } g \text{ is $\eta$-quasi-symmetric}\}.\)

It follows from the Arzela-Ascoli Theorem that $\widetilde{\mathcal F_\eta}$ is compact, so finitely many bi-Lipschitz approximations are sufficient to approximate any map in $\widetilde{\mathcal F_\eta}$ on $Q_0$. Given another tile $Q$, there are isometries $g,h$ of $H^+$ with $gQ=Q_0$, and the composition $gfh$ is in $\widetilde{\mathcal F_\eta}$.

Thus, in the first step, the PL approximation theorem is invoked only finitely many times. Furthermore, this produces only finitely many configurations (up to isometry) to interpolate in the second step. Thus, for each of the finitely many colors Sullivan's extension theorem is invoked only finitely many times.

The resulting map $F$ is then defined using finitely many bi-Lipschitz maps, and is therefore bi-Lipschitz and furthermore quasi-conformal. Since it is a bounded distance from $\widetilde  f$, $F$ has the prescribed extension to the boundary. Furthermore, the hyperbolic metric on $H^+$ is conformally equivalent to the Euclidean metric, so $F$ is a quasi-conformal map with respect to the Euclidean metric. Extending to all of $\R^{n+1}$ via reflection completes the proof.
\end{proof}
\end{thm}

Our goal now is to generalize Theorem \ref{thm:tv} as Theorem \ref{thm:mainmain}, applicable to a wider class of spaces, which we next define and study.

\subsection{Metric Similarity Spaces: Definition and Structure}
\label{sec:mss}
We first define stacked tilings. A metric space $H$ homeomorphic to $\R^n$ is said to have a \textit{stacked tiling} if it possesses both a discrete co-compact action by some group $\Gamma$ and a compatible homothety $\alpha$. More precisely:
\begin{defi}
\label{defi:stackedtiling}
Suppose $(H,d_H)$ is a complete metric space homeomorphic to $\R^n$ for some $n$ so that:
\begin{enumerate}
\item There exists a discrete group $\Gamma \subset \text{Isom}(H)$ acting properly on $H$ and a compact connected domain $K \subset H$ such that $H = \cup \{ \gamma K \st \gamma \in \Gamma\}$.
\item For all $\gamma \in \Gamma$ not equal to the identity, the interiors of $K$ and $\gamma K$ are disjoint.
\item There exists a homothety $\alpha: H \rightarrow H$ that rescales distances by some $a>1$.
\item \label{conjugationCondition}For $\gamma \in \Gamma$, $\alpha \gamma \alpha^{-1} \in \Gamma$.
\item \label{finitenessCondition}There is a finite subset $\Gamma' \subset \Gamma$ with $\alpha K = \cup \{\gamma' K \st \gamma' \in \Gamma'\}$.
\end{enumerate}
We then say that $H$ possesses a \textit{stacked tiling}, following the terminology of \cite{strichartz92}.
\end{defi}

\begin{example}
\label{ex:stackedtiling}
Consider $H=\R^n$ with group $\Gamma=\Z^n$. We may then take $K$ to be the unit cube and $\alpha(x) = 2x$. Then, $2K$ is a cube with side length 2, which is the union of $2^n$ translates of $K$.
\end{example}

\begin{question}
Is Condition \ref{finitenessCondition} equivalent to the subgroup $\alpha \Gamma \alpha^{-1}$ having finite index in $\Gamma$? 
\end{question}

\begin{defi}
\label{defi:metricsimilarityspace}
Suppose $H$ possesses a stacked tiling and furthermore $H^+ = H\times \R^+$ has a Riemannian metric $d_{H^+}$ so that
\begin{enumerate}
\item For $\gamma \in \Gamma$, the mapping $(x,t) \mapsto (\gamma x, t)$ of $H^+$ is an isometry.
\item The mapping $(x,t) \mapsto (\alpha x, a t)$ of $H^+$ is an isometry.
\end{enumerate}
Under these conditions, we refer to $H^+$ as an \textit{$(n+1)$-dimensional metric similarity space with base $H$}.
\end{defi}

\begin{example}
Let $H$ be as in Example \ref{ex:stackedtiling} and $H^+ = \R^n \times \R^+$. The hyperbolic metric on $H^+$ is then invariant under both the ``parabolic'' transformations $(x,t) \mapsto (x+x', t)$ and the ``hyperbolic'' transformation $(x,t) \mapsto (2x, 2t)$.
\end{example}

\begin{defi} 
Denote by $\Gamma_\alpha$ the semigroup
\[\Gamma_\alpha = \{ \alpha^n \gamma \st n \in \Z, \gamma \in \Gamma\}.\]
Depending on context, we view elements of $\Gamma_\alpha$ as mappings of either $H$ or $H^+$, based on Definition \ref{defi:metricsimilarityspace}. Condition \ref{conjugationCondition} of Definition \ref{defi:stackedtiling} provides a relation on the semigroup.
\end{defi}

If $\alpha$ commutes with elements of $\Gamma$, then $\Gamma_\alpha$ is a group. Even when this is the case, $\Gamma_\alpha$ is not discrete. However, the stacked tiling of $H$ can be turned into a tiling of $H^+$ by disjoint tiles related by $\Gamma_\alpha$. More specifically:

\begin{defi}
\label{defi:tiling}
Let $Q_0 = K \times [1, a]$. The \textit{tiles} of $H^+$ are the compact sets
\[\mathcal Q = \{ g Q_0 \st g \in \Gamma_\alpha\}.\]
We refer to $Q_0$ as the \textit{fundamental tile} of $H^+$. For each $n \in \Z$, $\{\alpha^n \gamma Q_0 \st \gamma \in \Gamma\}$ is the $n^{th}$ \textit{layer} of $\mathcal Q$. The tiles are isometric, and we have $\cup \mathcal Q = H^+$. Furthermore, the interiors of any two tiles are disjoint.  

A \textit{coloring} of $\mathcal Q$ is a partition $\mathcal Q = \mathcal Q_1 \sqcup \cdots \sqcup \mathcal Q_{N}$
such that two tiles of the same color (index) do not intersect. Because the tiles $\mathcal Q$ are in bijective correspondence with $\Gamma_\alpha$, every coloring of the tiles induces a coloring $\Gamma_1 \sqcup \cdots \sqcup \Gamma_{N}$ of $\Gamma_\alpha$.
\end{defi}

\begin{remark}
The fundamental tile $Q_0$ serves for a metric similarity space the role played by a fundamental domain in a space with a co-compact group action. However, even if $\Gamma_\alpha$ is a group, it is not discrete and $Q_0$ is not a fundamental domain.
\end{remark}

In general, elements of $\Gamma_\alpha$ do not take tiles to tiles.

\begin{defi}
Let $Q \in \mathcal Q$. We call an element $\gamma \in \Gamma_\alpha$ a \textit{combinatorial isometry with respect to $Q$} if for any tile $Q'$ with $Q \cap Q' \neq \emptyset$, the image $\gamma Q'$ of $Q'$ is a tile. We will not mention $Q$ explicitly if it is obvious which tile is relevant.
\end{defi}

\begin{lemma}
\label{lemma:layers}
Let $\gamma \in \Gamma$ and $Q \in \mathcal Q$. If $Q$ is in the $n^{th}$ layer of $\mathcal Q$, for $n\leq 0$, then $\gamma Q \in \mathcal Q$ is a tile.
\begin{proof}
A tile in the $n^{th}$ layer has the form $\alpha^n \gamma Q_0$. If $n \leq 0$, we can find using Condition \ref{conjugationCondition} of Definition \ref{defi:stackedtiling} $\gamma_1 \in \Gamma$ so that $\gamma \alpha^n \gamma Q_0 = \alpha^n \gamma_1 Q_0$, again a tile.
\end{proof}
\end{lemma}

\begin{lemma}[Normalization Lemma]
\label{lemma:finiteness}
Let $Q \in \mathcal Q$. Then $Q = g \gamma' Q_0$ for some $g \in \Gamma_\alpha$ and $\gamma' \in \Gamma'$, with $g$ a combinatorial isometry.
\begin{proof}

Because $\alpha$ takes tiles to tiles, we may assume $Q$ is in the $0^{th}$ layer of the tiling, so that $Q = \gamma_0 K \times [1, a]$ for some $\gamma_0 \in \Gamma$.

We now consider the stacked tiling of $H$. Because the tiles 
$\{ \gamma K \st \gamma \in \Gamma\}$ 
tile $H$, so do the tiles 
$\{ \alpha \gamma K \st \gamma \in \Gamma\}$. 
There then exists $\alpha \gamma_1 K$ with $\text{int}(\alpha \gamma_1 K) \cap \text{int}(\gamma_0 K) \neq \emptyset$, where int denotes the interior of a tile. Taking $g=\alpha \gamma_1 \alpha^{-1} \in \Gamma$, we have:
\(\alpha \gamma_1 K = g \alpha K = g \left(\bigcup_{\gamma' \in \Gamma'} \gamma'K\right).\)
Since $\text{int}(\alpha \gamma_1 K) \cap \text{int}(\gamma_0 K) \neq \emptyset$, we have $\gamma_0 = g \gamma'$.

Returning to $H^+$, we can write $Q= g \gamma' Q_0$. Tiles intersecting $Q$ are in either the same layer or one of the adjacent layers. By Lemma \ref{lemma:layers}, $g = \alpha \gamma_1 \alpha^{-1}$ takes tiles each of these three layers to other tiles in the same layer.
\end{proof}
\end{lemma}

\begin{cor}
\label{lemma:coloring}
There exists a coloring of $\mathcal Q$.
\begin{proof}
Every $Q \in \mathcal Q$ is equivalent to $\gamma' Q_0$ for some $\gamma' \in \Gamma'$ by a combinatorial isometry. There are finitely many tiles $\gamma'Q_0$, each with finite valence.
Thus, the valence of the adjacency graph is bounded, and it is therefore clear that it has a coloring.
\end{proof}
\end{cor}

\begin{cor}
\label{cor:disjointneighborhood}
Let $\Gamma_\alpha = \Gamma_1 \sqcup \cdots \sqcup \Gamma_N$ be a coloring induced from a coloring of $\mathcal Q$. Then there exists an open neighborhood $Q_1$ of $Q_0$ such that for any $i$ and distinct $g,h \in \Gamma_i$, $gQ_1 \cap h Q_1 = \emptyset$.
\begin{proof}
Let $\gamma' \in \Gamma'$. Let $D_{\gamma'} = d_{H^+}(\gamma' Q_0, (\gamma' Q_0)^*)$, where
\((\gamma' Q_0)^* = \cup \{ Q \in \mathcal Q \st Q \cap \gamma' Q_0 \neq \emptyset \}.\) Setting $D= \min_{\gamma' \in \Gamma'} D_{\gamma'} \neq 0$, let $Q_1$ be the $D/2$-neighborhood of $Q_0$. By Lemma \ref{lemma:finiteness}, we have considered all the tiles up to combinatorial isometry, so $Q_1$ has the desired property.
\end{proof}
\end{cor}

\subsection{Base and Boundary}
We next define a metric on $\mathcal Q$ using its adjacency graph and use a variation of the Milnor-\v Svarc Lemma (\cite{bridsonhaefliger} p. 140) to find a quasi-isometry between $H^+$ and $\text{Con}(H)$. This allows us to conclude that $H^+$ is hyperbolic and identify $H \cup \{\infty\}$ with $\partial_\infty H^+$.
\begin{defi}
Define a graph with vertices $\mathcal Q$ and an edge between two vertices $Q \neq Q'$ if $Q \cap Q' \neq \emptyset$. By an abuse of notation, we continue to refer to this graph as $\mathcal Q$. Denote by $d_{\mathcal Q}$ the graph metric on $\mathcal Q$ which assigns to each edge length $1$.
\end{defi}

\begin{lemma}
\label{lemma:structure}
Let $H^+=H \times \R^+$ be a metric similarity space and $\text{Con}(H)$ as in Definition \ref{defi:conx}. Then $id: H^+\rightarrow \text{Con}(H)$ is a quasi-isometry.
\label{lemma:milnor}
\begin{proof}
Recall that $d_{\text{Con}(H)}$ is $\epsilon$-roughly geodesic for some $\epsilon>0$. The metric $d_{H^+}$ is complete and Riemannian, so geodesic and in particular $\epsilon$-roughly geodesic. Let $d$ be either one of these metrics on the set $H \times \R^+ = H^+ = \text{Con}(H)$. Note that $\Gamma_\alpha$ acts on $(H^+,d)$ with isometries for either choice of $d$.  Increasing the fundamental tile $Q_0$ and taking the corresponding sub-semi-group of $\Gamma_\alpha$ if necessary, we may assume that $\epsilon \ll D = 3 \text{ diam}(Q_0)$. 

Define a function $f:H^+ \rightarrow \mathcal Q$ as follows. For each $x\in H^+$, let $f(x)\in \mathcal Q$ be a tile containing $x$. We now show that $f$ is a quasi-isometry between $(H^+, d)$ and the graph $\mathcal Q$. From this, we will be able to conclude that $id: (H^+, d_{H^+}) \rightarrow (\text{Con}(H), d_{\text{Con}(H)})$ is a quasi-isometry.

By definition of $D$, we know that if $x,y \in X$ and $d(f(x),f(y)) = 1$, then $d(x,y) \leq 2D/3$. Conversely, it is easy to see that there is an $N>0$ so that if $d(x,y)<D$, then $d_{\mathcal Q}(f(x),f(y)) \leq N$.

Let $x, y \in H^+$ and let $n = d_{\mathcal Q}(f(x), f(y))$. Pick a geodesic $f(x)=Q_1, \ldots, Q_{n+1} = f(y)$ in $\mathcal Q$, and choose points $q_1, \ldots, q_{n+1}$ in $H^+$ with $f(q_i)=Q_i$. We then have $d(q_i, q_{i+1}) \leq 2 D/3$ so by the triangle inequality we have
\(d(x,y)) \leq 2 n D/3 = 2 d_{\mathcal Q}(f(x),f(y)) D/3.\)
Thus, $f$ is co-Lipschitz.

Suppose now $x,y \in H^+$ are joined by an $\epsilon$-rough geodesic $r(s)$, with $s \in [0,d(x,y)]$. Let $r(s_0), \ldots, r(s_n)$ be a minimal sequence of points such that $s_0=0, s_n = d(x,y)$, and $d(r(s_i), r(s_{i+1})) \leq D$. We estimate
\(
d_{\mathcal Q}( f(x), f(y)) &\leq d_{\mathcal Q}(f(r(s_0)), f(r(s_1))) + \ldots + d_{\mathcal Q}(f(r(s_{n-1})), f(r(s_n)))\\
&\leq N n \leq N \left(\frac{d(x,y)}{D-\epsilon}+1\right),
\)
where the last inequality follows from the fact that $r$ is an $\epsilon$-rough geodesic and the choice of $s_i$ is minimal.

Since $f$ is onto, $\mathcal Q$ is quasi-isometric to both $d_{H^+}$ and $d_{\text{Con}(H)}$. Thus, the identity map from $H^+$ to $\text{Con(H)}$ is a quasi-isometry.
\end{proof}
\end{lemma}

\begin{lemma}
Let $f:X \rightarrow Y$ be a quasi-isometry between hyperbolic spaces. Suppose $d$ is in the conformal gauge of $\partial_\infty X$. Then the induced metric $(\partial_\infty f)_* d$ on $\partial Y$ is in the conformal gauge of $\partial Y$. The same holds for parabolic conformal gauges.
\label{lemma:conformalgauge}
\begin{proof}
Since $d$ is in the conformal gauge of $\partial_\infty X$, it is bi-Lipschitz to a visual metric $d_1$. Let $d_2$ be a visual metric on $\partial_\infty Y$. By Theorem \ref{thm:qiqs}, the induced metric $(\partial_\infty f)_* d_1$ on $\partial_\infty Y$ is quasi-symmetric to the metric $d_2$. Thus, same holds for $(\partial_\infty f)_*d$.
\end{proof}
\end{lemma}

The following follows directly from Lemma \ref{lemma:structure} and Lemma \ref{lemma:conformalgauge}:
\begin{thm}
\label{thm:boundary}
Let $H^+$ be a metric similarity space with base $H$. Then $H^+$ is Gromov hyperbolic and the mapping $\iota: H \cup \{\infty\} \rightarrow \partial_\infty H^+ $ defined in (\ref{fla:bdbijection}) is a bijection. Under this identification, the metric $d_H$ is in the parabolic conformal gauge of $\partial_\infty H^+\backslash \{[\gamma_\infty]\}$.
\label{thm:msshyperbolic}
\end{thm}

\begin{cor}
\label{cor:qiExtension}
Let $H^+ = H \times (0,\infty)$ be a metric similarity space. Under the identification $H \times \{0\} = \partial_\infty H^+ \backslash\{\infty\}$, every quasi-isometry of $H^+$ to itself extends continuously to a map on $H \times [0, \infty)$.
\end{cor}

\begin{remark}
\label{rmk:boundary}
In view of Theorem \ref{thm:msshyperbolic} and Corollary \ref{cor:qiExtension}, we use the notation $\partial H^+$ to denote both $\partial_\infty H^+$ and $H \cup \{\infty\}$. Likewise, for a quasi-isometry $f$ of $H^+$, $\partial f$ denotes the extension of $f$ to $\partial H^+$.
\end{remark}

\subsection{Extension from the Base}
\label{sec:lifts}

\begin{defi}
Let $H^+$ be an $n$-dimensional metric similarity space with base $H$.
A map $\wedge: C(H,H) \rightarrow C(H^+, H^+)$, denoted by $f \mapsto \widehat f$, is a \textit{lifting method} if
\begin{enumerate}
\item If $f$ is quasi-symmetric, then $\widehat f$ extends to a homeomorphism $\widehat f: H \times [0, \infty)\rightarrow H \times [0, \infty)$, with $\widehat f\vert_{H \times \{0\}} = f$.
\item $\wedge$ is continuous with respect to the compact-open topologies on the two spaces,
\item $\wedge$ is $\Gamma_\alpha$-equivariant: $\widehat{g_1 f g_2} = g_1 \widehat f g_2$ for all $g_1, g_2 \in \Gamma_\alpha$.
\end{enumerate}
\end{defi}

\begin{lemma}
\label{lemma:finitedistance}
Let $F$ be a quasi-isometry of $H^+$ preserving the vertical direction. If $\widehat{\partial F}$ is a quasi-isometry, then it is bounded distance from $F$.
\begin{proof}
Quasi-isometries are injective on the large scale, so there is a quasi-isometry $F^{-1}$ so that $F^{-1}\circ F$ is bounded distance from the identity map. Since $\partial (\widehat{\partial F}) = \partial F$, the boundary map $\partial (F^{-1} \circ \widehat{\partial F})$ is the identity map. It then follows from the stability of quasi-geodesics that $F^{-1} \circ \widehat{\partial F}$ is finite distance from the identity map.
\end{proof}
\end{lemma}

We now discuss a specific lifting method $\sim$ introduced by Tukia-\vaisala{}. Others are also available, as mentioned in the introduction.

\begin{defi}
Given $f: H \rightarrow H$ quasi-symmetric, let $\widetilde  f: H^+ \rightarrow H^+$ defined by 
\[\label{fla:basicextension}
&\widetilde f(x, t) = (f(x), \tau(x,t)) & \tau(x,t)=\max_{ \norm{x-y} \leq t}\norm{f(x)-f(y)}.
\]
\end{defi}

The following is easy to see from the definition of $\tau$:
\begin{lemma}
\label{lemma:homeo}The mapping $\sim: f \mapsto \widetilde f$ is a lifting method.
\end{lemma}

We next give a squence of examples demonstrating properties of $\sim$.
\begin{example}Let $f: \R^n \rightarrow \R^n$ be given by $f(x) = \lambda x$ for $\lambda \in \R$. Then $\widetilde f(x,t) = (\lambda x, \lambda t)$, a hyperbolic transformation of the upper half-space $\Hyp^{n+1}_\R$. More generally, if $F: \Hyp^{n+1}_\R \rightarrow \Hyp^{n+1}_\R$ is an isometry preserving the point at infinity, then $\widetilde{\partial F} = F$.
\end{example}

\begin{example}
Let $f: \R^2 \rightarrow \R^2$ be given by $f(x,y) = (x,2y)$. Then $\widetilde f(x,y,t) = (x,2y,2t)$ is 2-bi-Lipschitz, as can be seen by applying $\widetilde f$ to a geodesic and computing the effect on its speed. By extension, any affine map $f:\R^n \rightarrow \R^n$ lifts to a $\norm{\Lambda/\lambda}$-bi-Lipschitz map $\widetilde  f$, where $\Lambda$ is the eigenvalue of maximum norm and $\lambda$ is the eigenvalue of minimum norm.
\end{example}

\begin{example}
\label{ex:nobl}
Let $f: \R \rightarrow \R$ be given by $f(x) = x^3$, a quasi-symmetric map. Then $\widetilde  f: \Hyp^2_\R \rightarrow \Hyp^2_\R$ is given by $\widetilde  f(x,t) = \left(x^3, (\norm{x}+t)^3-\norm{x}^3\right)$. Away from $x=0$, $\widetilde f$ is differentiable, and the differential becomes singular as $x$ approaches $0$. Since the metric on $\Hyp^2_\R$ is Riemannian, this implies $\widetilde f$ is not bi-Lipschitz.
\end{example}

Example \ref{ex:nobl} shows that a quasi-symmetry $f$ does not necessarily lift to a bi-Lipschitz map under $\sim$. Theorem \ref{thm:main} implies, however, that for every quasi-symmetry $f$, the lift $\widetilde  f$ is a quasi-isometry and may in fact be approximated by a bi-Lipschitz map.

\subsection{Compactness of Quasi-Symmetries}
\label{sec:compact}
Let $H^+=H \times \R^+$ be a metric similarity space, and the notation as in Section \ref{sec:mss}. Fix a homeomorphism $\eta: [0,\infty) \rightarrow [0, \infty)$, distinct points $x_0, y_0 \in H$, and a lifting method $\wedge$.

\begin{defi}
A mapping $f: H \rightarrow H$ is \textit{normalized} if $f(x_0) \in K$ and 
\(
d_H(x_0, y_0) \leq d_H( fx_0,fy_0) \leq a d_H(x_0,y_0)
\)
for $a, K$ as in Definition \ref{defi:stackedtiling} for $H$.  We also define:
\( 
&\mathcal F_\eta = \left\{ f: H \rightarrow H \st f \text{ is a normalized $\eta$-quasi-symmetry} \right \} \subset C(H,H),\\
&\widehat{\mathcal F_\eta} = \left\{ \widehat f \st f \in \mathcal F_\eta\right\} \subset C(H^+,H^+).
\)
Note that $\widehat{\mathcal F_\eta}$ is the image of $\mathcal F_\eta$ under the continuous lifting map $\wedge$.
\end{defi}

\begin{remark}
\label{remark:normalization}
For any $\eta$-quasi-symmetric map $f: H\rightarrow H$, there exist $n \in \Z$ and $\gamma \in \Gamma$ so that $(\gamma \alpha^n f) \in \mathcal F_\eta$ and $(\gamma \alpha^n \widehat f) \in \widehat {\mathcal F_\eta}$. Note that $\gamma \alpha^n$ need not be in $\Gamma_\alpha$ if $n>0$.
\end{remark}

\begin{lemma}\label{compactness} $\mathcal F_\eta$ and $\widehat{\mathcal F_\eta}$ are compact in the compact-open topology.
\begin{proof}
It suffices to prove that $\mathcal F_\eta$ is compact, because $\widehat{\mathcal F_\eta}$ is its image under a continuous map. By the Arzela-Ascoli theorem for metric spaces, \cite{heinonen} Theorem 10.28, it suffices to show that $F_\eta$ is closed, equicontinuous and has compact orbits. We switch to the notation $\norm{x-y} := d_H(x,y)$ and follow the discussion on p. 85 of \cite{heinonen}.

The family $\mathcal F_\eta$ is defined by closed conditions, except possibly the condition that the maps be $\eta$-quasi-symmetric embeddings. It follows from the Arzela-Ascoli theorem that every limit of $\eta$-quasi-symmetric maps is a  $\eta$-quasi-symmetric embedding or constant. The condition $\norm{x_0-y_0} \leq \norm{f x_0 - f y_0}$ rules out the possiblity of a limit map that is constant. Thus $\mathcal F_\eta$ is closed.

For any $z \in H$ and $f\in \mathcal F_\eta$, we have by definition of quasi-symmetry:
\( \frac{\norm{fz - fx_0}}{\norm{fy_0-fx_0}} \leq \eta \left( \frac{\norm{z-x_0}}{\norm{y_0-x_0}}\right), \)
so that
\(\norm{fz-fx_0} \leq  \eta \left( \frac{\norm{z-x_0}}{\norm{y_0-x_0}}\right) \norm{fy_0-fx_0} \leq  a  \norm{x_0-y_0} \eta \left( \frac{\norm{z-x_0}}{\norm{y_0-x_0}}\right). \)
Since $f(x_0) \in K$, $\norm{fz-fx_0}$ is bounded for each $z \in H$. Since $H$ is complete and homeomorphic to $\R^n$, this implies that the closure of $\mathcal F_\eta (z)$ is compact for every $z\in H$.

We now prove equicontinuity at every $z_0 \neq x_0 \in H$. Let $z \in H$ and $f \in \mathcal F_\eta$.
\(\norm{ f z - f z_0} &\leq \eta \left(\frac{\norm{z-z_0}}{\norm{x_0-z_0}} \right) \norm{fx_0 - f z_0} \\
&\leq \eta \left(\frac{\norm{z-z_0}}{\norm{x_0-z_0}} \right)\eta \left(\frac{\norm{z-z_0}}{\norm{x_0-y_0}} \right) \norm{x_0-y_0}\\
&\leq \eta \left(\frac{\norm{z-z_0}}{\norm{x_0-z_0}} \right)\eta \left(\frac{\norm{x_0-z_0}}{\norm{x_0-y_0}} \right) \norm{f x_0-f y_0}
\\
& \leq a\norm{x_0-y_0} \eta \left(\frac{\norm{z-z_0}}{\norm{x_0-z_0}} \right)\eta \left(\frac{\norm{x_0-z_0}}{\norm{x_0-y_0}} \right).
\)
Thus, $\norm{fz-fz_0}$ is controlled by a  continuous function of $\norm{z-z_0}$. For $z_0 = x_0$, the same calculation holds with $x_0$ and $y_0$ switched. Thus, $\mathcal F_\eta$ is equicontinuous.
\end{proof}
\end{lemma}

\begin{cor}
\label{cor:finiteQS}
Up to specified additive error, $\widehat{\mathcal F_\eta}$ contains finitely many elements.
\end{cor}

\subsection{Bi-Lipschitz approximation}
\label{sec:bl}
We first re-phrase a standard result about Riemannian manifolds.

\begin{lemma}
\label{lemma:locallyEuclidean}
Let $H^+$ be a metric similarity space of dimension $n+1$. Then $H^+$ is locally bi-Lipschitz to $\R^{n+1}$.
\begin{proof}
Recall that $H^+$ is isometric to $\R^{n+1}$ with a Riemannian metric $d_{H^+}$. We claim that on any compact $K \subset \R^{n+1}$, $d_{H^+}$ is bi-Lipschitz to the Euclidean metric $d_E$. 

Let $D = \max \{ \text{diam}_{H^+}(K), \text{diam}_E(K) \}$, and let $K' = N_{H^+}(K, D) \cup N_E(K, D)$ denote the union of the closed $D$-neighborhoods of $K$ with respect to each metric. Note that $K'$ contains the union of all geodesic segments from $K$ to itself with respect to both metrics. Since both metrics are complete, $K'$ is compact. Next, denote by $T^1(K') \simeq K' \times S^{n-1}$ the (compact) unit tangent bundle to $K'$ as a subset of $\R^{n}$, with respect to either norm. Denote by $ds_E$ and $ds_{H^+}$ the line elements for the corresponding metrics and set
\(
L = \max \left\{ \sup_{v\in T^1(K')} \frac{ds_{H^+}(v)}{ds_E(v)}, \;\; \sup_{v\in T^1(K')} \frac{ds_E(v)}{ds_{H^+}(v)}\right\}.
\)
By definition of $d_E$ and $d_{H^+}$ as Riemannian metrics, they are $L$-bi-Lipschitz equivalent.
\end{proof}
\end{lemma}

The following follows immediately from standard piecewise-linear approximation results (see \cite{plapprox} p. 194) and Lemma \ref{lemma:locallyEuclidean}.
\begin{lemma}
\label{lemma:generalBLapprox}
Let $H^+$ be a metric similarity space of dimension $n \neq 4$, $f: {H^+} \rightarrow {H^+}$ a homeomorphism, and $K \subset {H^+}$ compact. Then for any $\epsilon > 0$ there exists a locally bi-Lipschitz $g: {H^+} \rightarrow {H^+}$ such that $\norm{g-f}_{K}<\epsilon$.
\end{lemma}

\begin{remark}
There exist counterexamples to Lemma \ref{lemma:generalBLapprox} in dimension 4.
\end{remark}

Lemma \ref{lemma:generalBLapprox} provides only local approximations of an embedding. Stitching local approximations together into a global bi-Lipschitz embedding is a difficult task. The only tool that is currently available for this purpose is the deep extension theorem of Sullivan. See \cite{sullivan79} for the original statement and \cite{tukiavaisala81} for a more detailed exposition of the proof. We state the Tukia-\vaisala{} adaptation of Sullivan's result in the context of metric similarity spaces.

\begin{lemma}[\cite{sullivan79, tukiavaisala82}]
\label{lemma:BL-extension} 
Let $H^+$ be a metric similarity space. Suppose $U \subset H^+$ has compact closure and $\mathcal F$ is a compact family of embeddings of $U$ into $H^+$ that have locally bi-Lipschitz approximations. Let $V \subset U$, $U' \subset U$ open with $\overline{U'} \subset U$, $V' \subset V$ with $\overline {V'} \subset V$. For any $\epsilon> 0$, there exists $\delta>0$ such that:

If $f \in \mathcal F$ and $g: V \rightarrow H^+$ is a bi-Lipschitz embedding with $d_{H^+}(f,g) < \delta$, then there is a bi-Lipschitz embedding $g': U \rightarrow H^+$ that extends the approximation of $f$ in the sense that:
\begin{enumerate}
\item $g'=g$ on $V'$,
\item $d_{H^+}(g', f) < \epsilon$ on $U'$.
\end{enumerate}
\end{lemma}

\begin{remark}
To paraphrase, if $g$ is an approximation of $f$ on $V$, then the approximation can be extended to $U$, at the cost of taking subsets of $U$ and $V$. The resulting error is at most $\epsilon$, as long as the original error was less than $\delta$. Note that $\delta$ depends on the compact family $\mathcal F$ in which $f$ is contained, but not on $f$ itself.
\end{remark}

\begin{remark}
For metric similarity spaces of dimension not equal to $4$, Lemma \ref{lemma:BL-extension} applies to the families $\widehat{\mathcal F_\eta}$, by Lemma \ref{compactness} and Lemma \ref{lemma:generalBLapprox}.
\end{remark}

Once we have approximated an injective map by a locally bi-Lipschitz map, we would like for it to remain an embedding. The following theorem guarantees that for every function sufficiently close to an element of $\widehat{\mathcal F_\eta}$, local injectivity implies global injectivity.

\begin{lemma}
\label{lemma:locallygloballyinjective}
Let $H^+$ be a metric similarity space, $r>0$, and $\widehat{\mathcal F_\eta}$ a family of normalized mappings of $H^+$. Then there exists $\epsilon > 0$ satisfying the following condition:\\
Let $h: {H^+} \rightarrow {H^+}$ with $\norm{h-\widehat f} < \epsilon$ for some $\widehat f \in \widehat{\mathcal F_\eta}$. Suppose that for any $z, w \in {H^+}$ with $\norm{z-w}<r$ we have $h(z) \neq h(w)$. Then $h$ is injective.
\begin{proof}
Recall that $Q_0$ denotes the compact fundamental tile of $H^+$ and that the family $\widehat{\mathcal F_\eta}$ is compact by Lemma \ref{compactness}. Since $H^+$ is a complete Riemannian space, its spheres $S(z,r) = \left\{w \in {H^+} \st d_{H^+}(z,w)\leq r\right\}$  are compact. We claim the lemma holds for
\(\epsilon = \frac{1}{3} \min \left\{ d_{H^+}\left(\widehat f z, \widehat f w\right) \st \widehat f \in \widehat{\mathcal F_\eta}, z\in Q_0, w \in S(z,r) \right\} > 0.
\)
Indeed, let $d_{H^+}\left(h, \widehat f\right) < \epsilon$, and $z,w \in H^+$. If $d_{H^+}(z,w)<r$, we are done, so we assume $d_{H^+}(z,w)\geq r$. If $z \in Q_0$, then for topological reasons we have 
\(
d_{H^+}\left(\widehat f z, \widehat f w\right) \geq d_{H^+}\left(\widehat f z, \widehat f \cdot S(z,r)\right) \geq 3 \epsilon.
\)
Since $d_{H^+}\left(\widehat f z, h z\right)< \epsilon$ and $d_{H^+}\left(\widehat f w, h w\right) < \epsilon$, we have $d_{H^+}(h(z),h(w)) \geq \epsilon$.

If $z \notin Q_0$, then $z \in g_1 Q_0$ for some $g_1 \in \Gamma_\alpha$. By Remark \ref{remark:normalization}, there is also an isometry $g_2$ of $H^+$ so that $g_2 f g_1^{-1} \in \widehat {\mathcal F_\eta}$. The calculation in the previous paragraph then applies, with $g_2 f g_1^{-1}$ substituted for $\widehat f$, $g_1 z$ for $z$, and $g_1 w$ for $w$.

We then have that for all $z,q\in H^+$ with $d_{H^+}(z,w)\geq r$, $d_{H^+}(hz, hw)>\epsilon$. In particular, $h$ is injective.
\end{proof}
\end{lemma}

\subsection{Main Theorem}
\label{sec:main}
We now prove the main theorem. For an overview of the proof in the case $H^+ = \Hyp^{n+1}_\R$, see Section \ref{sec:tv}.
\begin{thm}[Main Theorem] \label{thm:mainmain}
Let $H^+ = H \times \R^+$ be a metric similarity space of dimension not equal to 4, $\wedge$ a lifting method, $\epsilon > 0$, and $f: H \rightarrow H$ a quasi-symmetry. There exists a bi-Lipschitz embedding $F:H^+ \rightarrow H^+$ with $d_{H^+}(\widehat f, F) < \epsilon$, and $\partial F = f$.
\begin{proof}
Recall from Section \ref{sec:mss} that $H^+ = \cup \mathcal Q$, a union of tiles with disjoint interiors. The tiles $\mathcal Q$ are translates by the semigroup $\Gamma_\alpha$ of the fundamental tile $Q_0$. By Lemma \ref{lemma:coloring}, the tiling $\mathcal Q$ of $H^+$ is equipped with a coloring $\mathcal Q = \mathcal Q_1 \sqcup \cdots \sqcup \mathcal Q_{N}$, inducing a coloring $\Gamma_\alpha = \Gamma_1 \sqcup \cdots \sqcup \Gamma_{N}$.

We will approximate $\widehat f$ one tile at a time, first by piecewise-linear approximation on $\mathcal Q_1$, and using bi-Lipschitz extension on the other colors $i = 2, \ldots, N$. While there are infinitely many tiles of each color, we will use Lemma \ref{lemma:finiteness} to argue that PL approximation and bi-Lipschitz extension only need to be invoked finitely many times per color.

Step 1 (Base case: PL approximation). 
We will first approximate $\widehat f$ on $\mathcal Q_1$ within error $\epsilon_1$, to be specified later. By Corollary \ref{cor:disjointneighborhood}, there exists a neighborhood $Q_1$ of $Q_0$ with compact closure whose translates $\Gamma_1 Q_1$ are disjoint. We now apply PL approximation (Lemma \ref{lemma:generalBLapprox}) on the disjoint domains $\Gamma_1 Q_1$. 

Let $\gamma_1 \in \Gamma_1$. By Remark \ref{remark:normalization}, there exists $n \in \N$, $\gamma \in \Gamma$ so that $h \circ \widehat f \circ \gamma_1^{-1} \in \widehat{\mathcal F_\eta}$, where $h \gamma \alpha^n$ is an isometry of $H^+$. The family $\widehat{\mathcal F_\eta}$ is compact, so by Corollary \ref{cor:finiteQS} contains only finitely many maps (up to error $\epsilon/2$. Let $f'$ be one of these maps, equivalent to $h \circ \widehat f \circ \gamma_1^{-1}$ up to an error of $\epsilon/2$. Approximating $f'$ on $Q_1$ within error $\epsilon/2$ using Lemma \ref{lemma:generalBLapprox}, we get a bi-Lipschitz approximation $F_{\gamma_1}$ of $h \circ \widehat f \circ \gamma_1^{-1}$ with error at most $\epsilon_1$. Returning to the original tile, $h^{-1} \circ F_{\gamma_1} \circ \gamma_1$ approximates $\widehat f$ on the tile $\gamma_1 Q_1$.

We repeat in the same way on the remaining domains $\Gamma_1 Q_1$. The first approximation $F_1$ of $\widehat f$ is now defined on $\cup \Gamma_1 Q_1$, and approximates $\widehat f$ with error at most $\epsilon_1$.

Step 2 (Inductive step: Sullivan extension). Let $i>1$ and suppose $\widehat f$ has been approximated by a map $F_{i-1}$ on $\cup_{j=1}^{i-1} \Gamma_j Q_{i-1}$, where $Q_{i-1}$ is a neighborhood of $Q_0$ with compact closure. Assume furthermore that $F_{i-1}$ is defined by finitely many maps. More precisely, assume there is a finite collection of bi-Lipschitz maps $\mathcal F_{i-1}$ such that if $\gamma \in \Gamma_1 \cup \ldots \cup \Gamma_{i-1}$ and $h$ is an isometry of $H^+$ such that $h \circ \widehat f \circ \gamma^{-1} \in h \circ \widehat f \circ \gamma^{-1}$, then $h \circ F_{i-1} \circ \gamma^{-1} \in \mathcal F_{i-1}$.

Let $Q_i$ be a neighborhood of $Q_0$ in ${H^+}$ with compact closure with $\overline{Q_i} \subset Q_{i-1}$. By the choice of $Q_1$, the domains $\Gamma_i Q_i$ are disjoint.

Let $\gamma \in \Gamma_i$. The map $\widehat f$ is defined on all of $\gamma Q_{i-1}$. The previous approximation $F_{i-1}$ is defined on $\gamma Q_{i-1} \cap \left(\cup_{j=1}^{i-1} \Gamma_j Q_{i-1}\right)$. Applying Lemma \ref{lemma:BL-extension} with 
\(&U=\gamma Q_{i-1}, &V = \gamma Q_{i-1} \cap \left(\cup_{j=1}^{i-1} \Gamma_j Q_{i-1}\right)\\
&U' = \gamma Q_i,  &V' = \gamma Q_i \cap \left(\cup_{j=1}^{i-1} \Gamma_j Q_i\right)\)
and $\epsilon=\epsilon_i$ we get a bi-Lipschitz embedding $F_\gamma$ defined on $\gamma Q_i$ and compatible with $F_{i-1}$ on $\cup_{j=1}^{i-1} \Gamma_j Q_i$. Note that Lemma \ref{lemma:BL-extension} applies because of Remark \ref{remark:normalization} and Lemma \ref{lemma:generalBLapprox}. Repeat for all $\gamma \in \Gamma_i$ to get a new approximation $F_i$ of $\widehat f$ on $\cup_{j=1}^{i} \Gamma_j Q_j$. 

By the Finiteness Lemma \ref{lemma:finiteness} and the inductive finiteness assumption, we need apply Lemma \ref{lemma:BL-extension} only finitely many times. Thus, the finiteness assumption is preserved.

Step 3 (Completing the Induction).
By induction, $F=F_{N}$ is defined on all of ${H^+}$. Up to isometries of ${H^+}$, it is piecewise defined by elements of a finite family $\mathcal F_{N}$ of bi-Lipschitz maps. Thus, it is uniformly locally bi-Lipschitz. Note that a locally bi-Lipschitz map need not be globally injective.

Step 4 (Approximation Error and Injectivity). There are two demands on the approximation error. The first is the maximal error $\epsilon$ given by the statement of the theorem. Secondly, we would like for $F$ to be injective. Let $r = d_{H^+}( \partial Q_0, \partial Q_N)>0$. If $x,y \in H^+$ and $d(x,y) < r$, then $x$ and $y$ are both contained in an enlarged tile $\gamma Q_N$ for some $\gamma \in \Gamma_\alpha$. Thus, we know the approximation $F$ is injective on the scale $r$. Lemma \ref{lemma:locallygloballyinjective} provides an $\epsilon_N$ such that $F$ is globally injective if $d_H(F, \widehat f) < \epsilon_N$. We may assume $\epsilon_N< \epsilon$.

Inductively, suppose $\epsilon_i$ is already defined. Lemma \ref{lemma:BL-extension}, when applied in the $i^{th}$ step of the inductive argument above, provides an approximation of error at most $\epsilon_i$ as long as the previous error is at most $\epsilon_{i-1} < \epsilon_{i}$, for some $\epsilon_i$. We require this degree of error in the previous step.

The PL approxmation base case automatically complies with error requirement $\epsilon_1$.

With these error bounds on the approxmiation process, $F$ is injective. An injective uniformly locally bi-Lipschitz mapping of a path metric space is globally bi-Lipschitz.

Step 5 (Boundary). The maps $F$ and $\widehat f$ differ by bounded additive noise. It then follows from Theorem \ref{thm:boundary} that $\partial F = f$. 
\end{proof}

\end{thm}

\begin{cor}
\label{cor:explicitextension}
Let $H^+$ be a metric similarity space of dimension $n \neq 4$ with base $H$ and $f: H \rightarrow H$ a quasi-symmetry. Then for any lifting method $\wedge$, the mapping $\widehat f: H^+ \rightarrow H^+$ is a quasi-isometry.
\end{cor}

\begin{cor}
Let $H^+$ be a metric similarity space of dimension $n \neq 4$. Then any quasi-isometry $f$ of $H^+$ preserving the vertical direction is bounded distance from a bi-Lipschitz map. If the action of Isom($H^+$) on $\partial H^+$ is transitive, then $f$ need not preserve the vertical direction.
\begin{proof}
By Theorem \ref{thm:boundary}, $\partial f$ is a quasi-symmetry of $H$. Applying Theorem \ref{thm:mainmain} with any lifting method (e.g. the lifting method $\sim$) gives a bi-Lipschitz mapping $F$ of $H^+$ with $\partial F = \partial f$. It follows from the stability of quasi-geodesics that the distance between $f$ and $F$ is bounded.

If Isom($H^+$) acts on $\partial H^+$ transitively, then one may compose with the appropriate isometry in order to preserve the vertical direction.
\end{proof}
\end{cor}

\section{Discussion}
\label{sec:discussion}
In this section, we first present some examples of metric similarity spaces and then the connection to past results and open problems.

\subsection{Homogeneous Negatively Curved Spaces}

Let $M$ be a homogeneous negatively curved manifold. Heintze showed in \cite{heintze74} that $M$ is diffeomorphic to a Lie group $N \rtimes_{\alpha_s} \R^+$, where $N$ is some nilpotent group modeled on $\R^n$ and $\R^+$ acts on $\R^n$ by linear transformations $\alpha_s$. The group $N$ possesses a left-invariant metric for which $\alpha_s$ is a homothety of factor $s$. We first provide an example for $N = \R^n$. See Sections \ref{sec:rankone} and \ref{sec:xie} for additional examples. 

\begin{example}
\label{ex:twisted}
Let $\R^n \rtimes \R^+$ be a Lie group modeled on $\R^n \times \R^+$ with group structure
\[(x,t) \cdot (x', t') = (x+e^{At}x', t t'),\]
where $A$ is a diagonal matrix with positive eigenvalues $\lambda_1, \ldots, \lambda_n$. Like any Lie group, $\R^n \rtimes_A \R^+$ has a natural left-invariant metric. Define a metric $d_A$ on $\R^n$ by
\[&\norm{(x_1, \ldots, x_n)}_A = \norm{x_1}^{1/\lambda_1} + \cdots + \norm{x_n}^{1/\lambda_n} &d_A(x,x') = \norm{x-x'}_A.\]
\end{example}

\begin{prop}
\label{prop:twisted}
Suppose $A$ is a diagonal matrix with positive integer eigenvalues $\lambda_1, \ldots, \lambda_n$. Then the Lie group $\R^n \rtimes_A \R^+$ with any left-invariant Riemannian metric is a metric similarity space with base $(\R^n, d_A)$, 
\begin{proof}
Let $\Gamma = \Z^n$ and $K = [0,1]^n$. The mapping \(\alpha(x_1, \ldots, x_n) = (2^{\lambda_1}x_1, \ldots, 2^{\lambda_n}x_n)\) of $\R^n$ dilates distances by a factor of 2. The image of $K$ under $\alpha$ is 
\[\alpha K = [0,2^{\lambda_1}] \times \ldots \times [0, 2^{\lambda_n}],\]
tiled by finitely many copies of $K$. Furthermore, $\alpha$ commutes with $\Gamma$. Thus, $(\R^n, d_A)$ has a stacked tiling. The metric on $\R^n \rtimes_A \R^+$  is by constrution invariant under the transformations $\Gamma_\alpha$, so $\R^n \rtimes_A \R^+$ is a metric similarity space.
\end{proof}
\end{prop}

\begin{defi}
\label{ex:graded}
Let $G$ be a Lie group modeled on $\R^{n_1} \times \ldots \times \R^{n_r}$, with elements $x \in G$ written as $x = x_{ij}$ with $1 \leq i \leq r, 1 \leq j \leq n_i$, and group structure of the form
\[(x \cdot x')_{ij} = x_{ij} + x_{ij}' + F_{ij}(x,x').\]
The group $G$ is a \textit{rational graded nilpotent group of step $r$} if each $F_{ij}$ is a polynomial in the variables $x_1, \ldots, x_{i-1}, x_1', \ldots, x_{i-1}'$ with rational coefficients and satisfies
\[F_{ij}(\delta_s x, \delta_s x') = s^i F_{ij}(x,x'),\]
where the transformation $\delta_s$ is given by \[(\delta_s x)_{ij} = s^i x_{ij}.\]

The \textit{Carnot-Carath\`eodory} metric on $N$ is defined as follows. Let $\mathcal D$ be the left-invariant distribution that at the origin corresponds to the $\R^{n_1}$ subspace: $\mathcal D_0 = \R^{n_1}$. Call a path $\gamma(t)$ \textit{horizontal} if the velocity $\gamma'(t)\in \mathcal D_{\gamma(t)}$ for almost all $t$. Give $\mathcal D$ a left-invariant inner product $g_{\mathcal D}$. The Carnot-Carath\`eodory metric $d_{CC}$ between two points is defined as the $g_{\mathcal D}$-length of the shortest horizontal curve between the points. By Chow's Theorem, $d_{CC}$ is well-defined. It is clear that $\delta_s$ is a homothety of $d_{CC}$ by factor $s$.
\end{defi}

\begin{example}Both $\R^n$ and the Heisenberg group $\Heis^n$ (see Section \ref{sec:rankone}) are rational graded nilpotent groups. See \cite{strichartz92} for further examples.
\end{example}

\begin{thm}[Strichartz \cite{strichartz92}] 
\label{thm:strichartz}Every rational graded nilpotent group possesses a stacked tiling.
\end{thm}

\begin{cor}
Let $N$ be a rational graded nilpotent group. Let $\R^+$ act on $N$ by the dilations $\delta_s$ in Definition \ref{ex:graded}. Then $N \rtimes_{\delta_s} \R^+$ with any left-invariant Riemannian metric is a metric similarity space.
\begin{proof}
The group $N \rtimes_{\delta_s} \R^+$ acts on itself by isometries. By \ref{thm:strichartz}, if one picks the homothety $\alpha = \delta_2$ of $N$, there is a corresponding subgroup $\Gamma \subset N$ and a compact $K \subset N$ satisfying the definition of stacked tiling.
\end{proof}
\end{cor}

\begin{remark}
Every homogeneous negatively curved manifold decomposes as $N \rtimes_{\alpha_s} \R^+$ for some family of diffeomorphisms $\alpha_s$. However, even if $N$ is a rational graded nilpotent group, generally $\alpha_s\neq \delta_s$. See Example \ref{ex:twisted}.
\end{remark}

\subsection{Rank One Symmetric Spaces}
\label{sec:rankone}
Among the homogeneous negatively curved manifolds, the most-studied ones are the non-compact rank one symmetric spaces. These are the real, complex, and quaternionic hyperbolic spaces, as well as the octonionic plane. In this section, we show that both real hyperbolic space $\Hyp^{n+1}_\R$ and complex hyperbolic space $\Hyp^{n+1}_\C$ are metric similarity spaces. The same is true of the quaternionic and octonionic spaces. However, according to Pansu rigidity \cite{pansu-metriques}, any quasi-isometry of the quaternionic and Cayley spaces is in fact an isometry, and any quasi-symmetry of the boundary is conformal. There exist quasi-isometries of $\Hyp^{n+1}_\R$ and $\Hyp^{n+1}_\C$ that are not equivalent to isometries (see \cite{Koranyi-Reimann1998}).

\begin{defi}
Real hyperbolic space $\Hyp^{n+1}_\R$ is the space $\R^n \rtimes_A \R^+$ defined in Example \ref{ex:twisted}, where $A$ is taken to be the identity matrix. The line element is explicitly given by 
\(ds^2 = \frac{dx^2 + dt^2}{t^2}.\)
\end{defi}

\begin{prop}
Real hyperbolic space $\Hyp^{n+1}_\R$ is a metric similarity space with base $(\R^n,d_E)$, where $d_E$ is the Euclidean metric.
\begin{proof}
This follows directly from Proposition \ref{prop:twisted}.
\end{proof}
\end{prop}

We now describe complex hyperbolic space $\Hyp^n_\C$. The definition corresponds to that of the Klein disk model of $\Hyp^n_\R$, but uses complex coefficients instead of real ones. For more details, see the book of Goldman \cite{goldman} and a recent survey of Parker-Platis \cite{parker-platis2010}.

Denote by $\CP^n$ complex projective space of dimension $n$. Let $J$ be a Hermitian form on $\C^{n+1}$ of type $(n,1)$, so $J(z,z) \in \R$ for all $z$. Complex hyperbolic space is the set of negative points in $\CP^n$ with respect to $J$ (we will define a metric shortly).
\(\Hyp^n_{\C,J} = \{p \in \CP^n \st J(z,z)<0 \text{ for } z\in \C^{n+1} \text{ representing } p\}.\)
Because $J$ is unique up to a change of basis, the space is well-defined up to a choice of coordinates, and we write $\Hyp^n_\C$ if $J$ is not specified. Taking $z = (z_0, \ldots, z_n)$, we define two Hermitian forms:
\( J_1(z) &= -\norm{z_0}^2 + \norm{z_1}^2 + \ldots + \norm{z_n}^2, \\
J_2(z) &= z_0 \overline{z_1} + z_1 \overline{z_0}+ \norm{z_2}^2 + \ldots + \norm{z_n}^2.
\)
In the coordinate patch $z_0 = 1$, $\Hyp^n_{\C,J_1}$ is the unit ball in $\C^n$, while $\Hyp^n_{\C,J_2}$ is the Siegel (or paraboloid) model of complex hyperbolic space:
\(
\Hyp^n_{\C,J_2} = \{ (z_1, \ldots, z_n) \st   -2 \Re(z_1)> \norm{z_2}^2 + \ldots + \norm{z_n}^2\}.
\)
Let $z,w \in \C^{n+1}$ represent two points $p,q \in \Hyp^n_{\C,J}$. The metric on $\Hyp^n_{\C,J}$ is defined by
\(
\cosh^2 d_{\Hyp^n_{\C,J}}(p,q) = \frac{J(z,w)J(w,z)}{J(z,z)J(w,w)}.
\)
This metric is, in fact, Riemannian. It may alternately be defined by taking the standard inner product at the origin of the ball model $\Hyp^n_{\C, J_1}$ and uniquely extending to a metric tensor on $\Hyp^n_{\C,J_1}$ invariant under the group of projective transformations preserving $J_1$.

\begin{remark}Unlike real hyperbolic space, complex hyperbolic space does not have constant curvature: for $n>1$ the sectional curvature $\kappa$ of $\Hyp^n_\C$ varies in the range $-4 \leq \kappa \leq -1$.  In the case $n=1$, $\Hyp^1_{\C,J_1}$ is precisely the Poincar\`e model of the real hyperbolic plane, rescaled to have constant curvature $-4$.
\end{remark}

Let $\Heis^{n-1}$ denote the 2-step nilpotent Heisenberg group (see also \cite{tysonetal, Koranyi-Reimann1995}) modeled on $\C^{n-1} \times \R$ with group law
\[(\zeta,u) \cdot (\zeta', u') = (\zeta+\zeta', u+u'+2\Im \langle \zeta, \zeta'\rangle),\]
where $\langle \zeta, \zeta'\rangle$ is the standard Hermitian pairing of complex vectors. The Heisenberg group is a rational graded nilpotent group of step 2. As such, it has a Carnot-Carath\`eodory metric with dilations $\delta_s$ (see Definiton \ref{ex:graded}).

Define the horospherical model of complex hyperbolic space as $\Hyp^n_{\C,\text{horo}} = \Heis^{n-1} \times \R^+ = \C^{n-1} \times \R \times \R^+$.
Give $\Hyp^n_{\C, \text{horo}}$ a Riemannian metric by identifying it with the Siegel model by means of the diffeomorphism
\(&f: \Hyp^n_{\C, \text{horo}} \rightarrow \Hyp^n_{\C,J_2} \subset \C \times \C^{n-1} &f:(\zeta,u,v) \mapsto \left(-v + \ii u + \norm{\zeta}^2, \sqrt{2}\zeta\right).\notag
\)
In horospherical coordinates, we have the following relationship between $\Hyp^n_\C$ and $\Heis^{n-1}$:
\begin{enumerate}
\item Isometries of $\Heis^{n-1}$ extend trivially to parabolic isometries of $\Hyp^n_{\C, horo}$.
\item Similarities $\delta_s$ of $\Heis^{n-1}$ extend to hyperbolic isometries defined by $(\zeta, u,v) \mapsto (s \zeta, s^2 u, s^2 v)$.
\end{enumerate}
Furthermore, we have the following properties of $\Heis^{n-1}$.
\begin{enumerate}
\item The subgroup $\Gamma\subset \Heis^n$ consisting of elements with integer coordinates is discrete and co-compact. 
\item By Theorem \ref{thm:strichartz} or the results of \cite{baloghtyson2006}, there exists a ``Heisenberg cube'', a fundamental domain for $\Gamma$ that provides $\Heis^{n-1}$ with a stacked tiling with $\alpha=\delta_2$.
\end{enumerate}

\begin{prop}
Complex hyperbolic space $\Hyp^n_\C$ is a metric similarity space with base $\Heis^{n-1}$ with its Caranot-Carath\`eodory metric.
\begin{proof}
The horospherical model $\Hyp^n_\C$ satisfies the definition, up to taking the square root of the $v$ coordinate, which does not scale properly as stated.
\end{proof}
\end{prop}

\subsection{Comparison with a Result of Xie}
\label{sec:xie}

A Hadamard manifold is a complete negatively curved manifold with pinched sectional curvature. A general Hadamard manifold is not a metric similarity space. However, the following result has significant overlap with our main theorem:

\begin{thm}[Xie \cite{xie2009}]
\label{thm:xie}
Let $M$ be a Hadamard manifold that is the universal cover of a compact manifold of dimension not equal to $4$. Then every quasi-isometry from $M$ to itself is bounded distance from a bi-Lipschitz map of $M$ to itself.
\end{thm}

In particular, Theorem \ref{thm:xie} applies to complex hyperbolic space and other negatively curved homogeneous manifolds admitting co-compact discrete group actions.
We present two examples to distinguish our main Theorem \ref{thm:mainmain} from Theorem \ref{thm:xie}.

The following well-known example shows that Theorem \ref{thm:xie} does not apply to most metric similarity spaces constructed in Example \ref{ex:twisted}:
\begin{prop}
\label{prop:badlytwisted}
Let $\R^n \rtimes_A \R^+$ be as in Proposition \ref{prop:twisted}, with the $\lambda_i$ distinct positive integers. Then $\R^n \rtimes_A \R^+$ does not admit a co-compact isometric group action.
\begin{proof}
Consider the group of isometries Isom$(\R^n \rtimes_A \R^+)$. It is easy to see by considering the Riemannian metric on $\R^n \rtimes_A \R^+$ that in this case the stabilizer Stab$(p)$ of a point is finite.

Suppose that a discrete group $G$ acts co-compactly on $\R^n \rtimes_A \R^+$. There would then be a finite-index subgroup $G'$ of $G$ whose elements are all either parabolic or hyperbolic. In particular, the left-invariant vector field along the $\R^+$ direction would descend to $\R^n \rtimes_A \R^+/G'$. The resulting flow on $\R^n \rtimes_A \R^+/G'$ would be isometric but not volume preserving, which is impossible.
\end{proof}
\end{prop}

The next proposition provides a perturbation method for metric similarity spaces.
\begin{prop}
\label{prop:bump}
Let $H^+$ be a metric similarity space with base $H$, a Riemannian metric $g$, and  $Q_0 = K \times [1,a]$. Let $f$ be a positive bump function on the interval $[1,a]$, equal to 1 near the endpoints. Extend $f$ to a function on $(0, \infty)$ invariant under the mapping $x \mapsto ax$. Then $(H^+, fg)$ is a metric similarity space. 
\begin{proof}
If $g$ is complete and invariant under the semigroup $\Gamma_\alpha$, then so is $fg$.
\end{proof}
\end{prop}

\begin{remark}
A more interesting version of Proposition \ref{prop:bump} would allow arbitrary perturbations of the Riemannian metric within the fundamental tile $Q_0$. For the sake of simplicity, we required $\Gamma_\alpha$ to act on $H^+$ by isometries when we  defined metric similarity spaces. However, our results remain valid if $g \in \Gamma_\alpha$ is only required to be an isometry on some fixed neighborhood of $Q_0$. Thus, one may introduce arbitrary bumps in the metric, including the appearance of positive curvature.
\end{remark}

\subsection{Concluding Remarks}

Theorem \ref{thm:mainmain} states that every quasi-symmetry of the base of a metric similarity space of dimension not equal to 4 extends to a bi-Lipschitz map of the space itself. This leaves open the following problem.
\begin{prob}
\label{prob:dim4}
Let $H^+$ be a 4-dimensional metric similarity space with base $H$. Is it true that every quasi-symmetry of $H^+$ preserving the vertical direction is bounded distance from a bi-Lipschitz map? Equivalently, does every quasi-symmetry of $H$ lift to a bi-Lipschitz map of $H^+$?
\end{prob}

The exclusion of 4-dimensional spaces from the main theorem is due to the dimension restriction in Lemma \ref{lemma:generalBLapprox}. Because there do exist homeomorphisms of $\R^4$ that cannot be approximated piecewise-linearly (see \cite{plapprox}), the restriction cannot be removed. However, such pathological homeomorphisms are difficult to construct. This suggests the following problem.

\begin{prob}
Let $H^+$ be a 4-dimensional metric similarity space with base $H$. Let $\wedge$ be a lifting method, $\epsilon>0$, and $f: H \rightarrow H$ a quasi-symmetric embedding. Can $\widehat f$ be approximated by a bi-Lipschitz map within error $\epsilon$?
\end{prob}

In the case of $H^+  = \Hyp^4_\R$, one bypasses the dimension restriction of Lemma \ref{lemma:generalBLapprox} using PL approximation in $\R^3$ and the following property of the lifting method $\sim$ defined by (\ref{fla:basicextension}).

\begin{lemma}[Tukia-\vaisala{} \cite{tukiavaisala82}]
\label{lemma:BLlift}
Let $f: \R^n \rightarrow \R^n$ be a homeomorphism. Then $f$ is locally bi-Lipschitz in the Euclidean metric if and only if $\widetilde f$ is locally bi-Lipschitz in the hyperbolic metric on $\Hyp^{n+1}_\R$.
\end{lemma}

Lemma \ref{lemma:BLlift} does not generalize directly to the setting of metric similarity spaces, as its proof uses the fact that the hyperbolic metric is locally bi-Lipschitz to the $l^1$ sum of the Euclidean metrics on $\R^n$ and $\R^+$. 

\begin{prob}
\label{prob:generalize}
Let $H^+$ be a metric similarity space with base $H^+$. Generalize Lemma \ref{lemma:BLlift} for some lifting method $\wedge$. More specifically, relate the analytic properties of $f: H \rightarrow H$ and those of $\widehat f: H^+ \rightarrow H^+$.
\end{prob}

Shanmugalingam-Xie showed in \cite{shanmugalingam-xie2011} that for many nilpotent groups $N$, all quasi-symmetries of $N$ are in fact bi-Lipschitz. While quasi-symmetries of the Heisenberg group need not be bi-Lipschitz, Capogna and Tang constructed a large class of bi-Lipschitz quasi-symmetries of $\Heis^n$ in \cite{capogna-tang}. An appropriate analogue of $\ref{lemma:BLlift}$ would allow these to be lifted to bi-Lipschitz mappings of $\Hyp^{n+1}_\C$. The existence of quasi-symmetries of $\Heis^n$ that are not bi-Lipschitz suggests the following problem.

\begin{prob}
Develop an analogue of piecewise-linear approximation theory for $\Heis^n$. Following the Pansu-Rademacher theorem \cite{pansu-metriques}, the phrase ``piecewise linear'' should generalize to ``piecewise horizontally homomorphic'', where a horizontal homomorphism of $\Heis^n$ preserves the standard splitting of the tangent spaces of $\Heis^n$.
\end{prob}

One may also attempt a more direct approach to proving the main theorem for $H^+ = \Hyp^{2}_\C$. 
Beurling-Ahlfors originally proved Theorem \ref{thm:tv} in the case of $H^+=\Hyp^2_\R$ \cite{beurling-ahlfors} using an explicit integral formula. While Theorem \ref{thm:tv} is proven using a different approach, perhaps an explicit bi-Lipschitz lifting may be found for $H^+=\Hyp^2_\C$.

Ahlfors' proof in \cite{ahlfors64} of Theorem \ref{thm:tv} in the case of $\Hyp^3_\R$ was based on a decomposition of a quasi-conformal mappings. Given a quasi-conformal mapping $f$ of $\R^2$, one finds a quasi-conformal flow $f_t$ on $\R^2$ so that $f = f_1$. One can then write $f$ as the composition of mappings with small distortion, which can then be approximated by bi-Lipschitz maps using results concerning almost-conformal mappings.

\begin{prob}
For $n\geq 3$, can a quasi-conformal mapping of $\R^n$ be embedded in a quasi-conformal flow? For $n \geq 1$, can a quasi-conformal mapping of $\Heis^n$ be embedded in a quasi-conformal flow? See \cite{Koranyi-Reimann1998} for a study of quasi-conformal flows on $\Heis^n$ and their lifts to $\Hyp^{n+1}_\C$.
\end{prob}

\bibliography{bib}

\end{document}